\documentclass[leqno]{amsart}

\usepackage[leqno]{amsmath}
\usepackage{amssymb,amsthm}
\usepackage{enumerate, latexsym, mathrsfs, mdwlist}

\swapnumbers
\theoremstyle{definition}
\newtheorem{nul}{}[section]
\newtheorem{dfn}[nul]{Definition}

\newtheorem*{dfn*}{Definition}
\newtheorem*{axm*}{Axiom}
\newtheorem*{ntn*}{Notation}
\newtheorem*{exm*}{Example}
\newtheorem*{exr*}{Exercise}
\newtheorem*{int*}{Intuition}
\newtheorem*{qst*}{Question}

\theoremstyle{plain}

\newtheorem{thm}[nul]{Theorem}
\newtheorem{prp}[nul]{Proposition}

\newtheorem{cor}[nul]{Corollary}
\newtheorem*{thm*}{Theorem}
\newtheorem*{prp*}{Proposition}
\newtheorem*{cor*}{Corollary}
\newtheorem*{lem*}{Lemma}
\newtheorem*{cnj*}{Conjecture}

\numberwithin{equation}{nul}

\newcommand{\PP}{\mathbf{P}}
\newcommand{\QQ}{\mathbf{Q}}

\renewcommand{\SS}{\mathbf{S}}

\newcommand{\ZZ}{\mathbf{Z}}

\newcommand{\tmf}{\mathrm{tmf}}
\newcommand{\TMF}{\mathrm{Tmf}}
\newcommand{\ko}{\mathrm{ko}}
\newcommand{\KO}{\mathrm{KO}}
\newcommand{\ku}{\mathrm{ku}}
\newcommand{\KU}{\mathrm{KU}}
\newcommand{\xxi}{\overline{\xi}}
\newcommand{\BP}[1]{\mathrm{BP}\langle{#1}\rangle}

\newcommand{\coloneq}{\mathrel{\mathop:}=}

\makeatletter 
\def\revddots{\mathinner{\mkern1mu\raise\p@ 
\vbox{\kern7\p@\hbox{.}}\mkern2mu 
\raise4\p@\hbox{.}\mkern2mu\raise7\p@\hbox{.}\mkern1mu}} 
\makeatother

\usepackage{tikz}
\usetikzlibrary{matrix,arrows,decorations}

\pgfdeclaredecoration{single line}{initial}{\state{initial}[width=\pgfdecoratedpathlength-1sp]{\pgfpathmoveto{\pgfpointorigin}}\state{final}{\pgfpathlineto{\pgfpointorigin}}} 

\newcommand{\fromto}[2]{{#1}\ \tikz[baseline]\draw[>=stealth,->](0,0.5ex)--(0.5,0.5ex);\ {#2}}

\newcommand{\goesto}[2]{{#1}\ \tikz[baseline]\draw[|->](0,0.5ex)--(0.5,0.5ex);\ {#2}}

\renewcommand{\to}{\ \tikz[baseline]\draw[>=stealth,->](0,0.5ex)--(0.5,0.5ex);\ }

\usepackage{color}

\title{Regularity of structured ring spectra and localization in $K$-theory}

\author{Clark Barwick}
\address{Massachusetts Institute of Technology, Department of Mathematics, Building 2, 77 Massachusetts Avenue, Cambridge, MA 02139-4307, USA}
\email{clarkbar@gmail.com}

\author{Tyler Lawson}
\address{University of Minnesota, School of Mathematics, Vincent Hall, 206 Church St. SE, Minneapolis, MN 55455, USA}
\email{tlawson@math.umn.edu}

\begin{document}

\maketitle

An early success of Quillen's algebraic $K$-theory was his localization sequence \cite[\S 5]{quillen}.  This sequence shows that algebraic $K$-theory enjoys a strong excision property that permits one to cut out certain subvarieties. For example, suppose $R$ a regular noetherian ring, and suppose $x\in R$ an element such that the quotient $R/x$ is also regular.  Then the localization sequence takes the form of a long exact sequence
\begin{equation*}
\cdots\to K_n(R/x)\to K_n(R)\to K_n(R[x^{-1}])\to K_{n-1}(R/x)\to\cdots.
\end{equation*}
The aim of this short paper is to identify a \emph{regularity} property for structured ring spectra (\textsc{aka} associative $\SS$-algebras, \textsc{aka} $E_1$ rings, \textsc{aka} $A_\infty$ rings) that allows one to prove a natural analogue of Quillen's theorem in this context.

To obtain Quillen's sequence, one may begin with a localization sequence
\begin{equation*}
\cdots\to K_n(\mathrm{Nil}_{(R,x)})\to G_n(R)\to G_n(R[x^{-1}])\to K_{n-1}(\mathrm{Nil}_{(R,x)})\to\cdots
\end{equation*}
where $\mathrm{Nil}_{(R,x)}$ is the category of finitely generated, $x$-nilpotent $R$-modules. This sequence is exact irrespective of any regularity hypotheses on $R$ and $x$; for example, it is an instance of the Fibration Theorem of Waldhausen \cite[Th. 1.6.4]{MR86m:18011}. The regularity condition actually enters twice to convert this sequence into the sequence above. First, one deduces isomorphisms
\begin{equation*}
G_{\ast}(R/x)\cong K_{\ast}(R/x)\textrm{,\quad}G_{\ast}(R)\cong K_{\ast}(R)\textrm{,\quad and\quad}G_{\ast}(R[x^{-1}])\cong K_{\ast}(R[x^{-1}])
\end{equation*}
from corresponding equivalences of derived categories. Second, one uses Quillen's D\'evissage Theorem \cite[Th. 4]{quillen} to identify $K_{\ast}(\mathrm{Nil}_{(R,x)})$ and $G_{\ast}(R/x)$.

What's remarkable about Quillen's D\'evissage Theorem is that it provides an equivalence of $K$-theories that is not induced by an equivalence of derived categories. Results of this kind are rare commodities, but the work of Blumberg and Mandell \cite{ku} provides another D\'evissage Theorem: they show that for any connective $E_1$ ring $\Lambda$, the algebraic $K$-theory of the category of $\Lambda$-modules with finitely generated homotopy and finite Postnikov towers is the same as the algebraic $K$-theory of $\pi_0\Lambda$. A version of the Blumberg--Mandell D\'evissage Theorem can, for example, be deduced from the Theorem of the Heart for Waldhausen $K$-theory \cite{heart}, which states
that the inclusion $\mathscr{A}^{\heartsuit}\subset\mathscr{A}$ of the
heart of a bounded t-structure on a stable $\infty$-category
$\mathscr{A}$ induces an isomorphism
\begin{equation*}
K_{\ast}(\mathscr{A}^{\heartsuit})\cong K_{\ast}(\mathscr{A}).
\end{equation*}

In this paper, we introduce a notion of \emph{regularity} on $E_1$ rings. Then, for a regular $E_1$ ring $\Lambda$, we use the Blumberg--Mandell D\'evissage Theorem of to identify the fiber term of a localization sequence
\begin{equation*}
K(\mathrm{Acyc}(L))\to K(\Lambda)\to K(L(\Lambda))
\end{equation*}
as the $K$-theory of the ordinary ring $\pi_0\Lambda$. Here $L$ is any smashing localization on $\Lambda$-modules with the property that a perfect $\Lambda$-module is $L$-acyclic if and only if it has only finitely many nontrivial homotopy groups. The result is a long exact sequence
\begin{equation*}
\cdots\to K_n(\pi_0\Lambda)\to K_n(\Lambda)\to K_n(L(\Lambda))\to K_{n-1}(\pi_0\Lambda)\to\cdots
\end{equation*}
The value of such a result is that it provides a description of the $K$-theory of a nonconnective ring spectrum $L(\Lambda)$ entirely in terms the $K$-theory of connective ring spectra.

Let us now briefly describe regularity for (suitably finite) connective structured ring spectra $\Lambda$. We'll say that a left $\Lambda$-module $M$ is \emph{coherent} if it has only finitely many nonzero homotopy groups, all of which are finitely generated as $\pi_0\Lambda$-modules. Now $\Lambda$ will be said to be \emph{regular} if $\pi_0\Lambda$ is a regular (commutative) ring, and every coherent left $\Lambda$-module is a retract of a finite cell module. If $\Lambda$ is regular, then the $G$-theory of $\Lambda$ --- which is the $K$-theory of the category of coherent $\Lambda$-modules --- maps naturally to the $K$-theory of $\Lambda$, and the D\'evissage Theorem of Blumberg--Mandell identifies $G_{\ast}(\Lambda)$ with $G_{\ast}(\pi_0\Lambda)$, which is in turn $K_{\ast}(\pi_0\Lambda)$. This is the key to identifying the fiber term of the fiber sequence above.

Our regularity condition sounds quite abstract, but in fact it can be checked entirely in terms of $\pi_0\Lambda$-modules: $\Lambda$ is regular just in case $\pi_0\Lambda$ is regular and $H\pi_0\Lambda\wedge_{\Lambda}H\pi_0\Lambda$ is a retract of a finite cell $H\pi_0\Lambda$-module. This condition allows us to check regularity in several examples.

First, for real $K$-theory, we have
\begin{equation*}
K(\ZZ)\to K(\ko)\to K(\KO).
\end{equation*}
Andrew Blumberg has informed us that Vigleik Angeltveit originally demonstrated how to prove this using the results of \cite{ku} at a conference in 2008 \cite{birs}. It confirms a conjecture of Christian Ausoni and John Rognes.

Next, for topological modular forms, we have a fiber sequence
\begin{equation*}
K(\ZZ)\to K(\tmf)\to K(\TMF).
\end{equation*}
Here $\tmf$ is the connective topological modular forms spectrum, and $\TMF$ is the nonconnective, nonperiodic version associated to the compactified moduli of elliptic curves.  We are unsure whether this fiber sequence was expected by others.

Finally, let $\BP{n}$ denote a truncated Brown--Peterson spectrum, and $\BP{n}^{\ast}$ denote the homotopy limit of the $E_1$ rings $S^{-1} \BP{n}$ for $S\subset\{v_1,\dots,v_n\}$, $S \neq \emptyset$.  (This is the localization away from $(v_1,\ldots,v_n)$ in the sense of \cite{gmcompletions}.) We obtain fiber sequences
\begin{equation*}
K(\ZZ_{(p)})\to K(\BP{n})\to K(\BP{n}^{\ast}).
\end{equation*}
In particular, when $n=1$ this specializes to a $p$-local version of the localization sequence of Blumberg-Mandell \cite{ku}.

\subsection*{Acknowledgements} We thank Andrew Blumberg and John Rognes for their thorough comments and corrections.


\section{Regularity} In effect, regularity for $E_1$ rings is a finiteness property. There is a host of finiteness properties for $E_{1}$ rings and their modules that one may consider. We quickly call a few of these to mind here. For more details, we refer to Jacob Lurie's discussion in \cite[\S 8.2.5]{HA}

Recall \cite[Pr. 8.2.5.16]{HA} that a connective $E_1$ ring $\Lambda$ is said to be \emph{left coherent} if $\pi_0\Lambda$ is left coherent as an ordinary ring, and if for any $n\geq 1$, the left $\pi_0\Lambda$-module $\pi_n\Lambda$ is finitely presented.

A left module $M$ over a left coherent $E_1$ ring $\Lambda$ is \emph{almost perfect} just in case $\pi_mM=0$ for $m\ll 0$ and for any $n$, the left $\pi_0\Lambda$-module $\pi_nM$ is finitely presented \cite[Pr. 8.2.5.17]{HA}.

\begin{dfn} Suppose $\Lambda$ a left coherent $E_1$ ring, and suppose $M$ a left $\Lambda$-module. We say that $M$ is \emph{truncated} if $\pi_mM=0$ for $m\gg 0$. We will say $M$ is \emph{coherent} if it is both truncated and almost perfect.
\end{dfn}

Recall \cite[Pr. 8.2.5.21]{HA} that a left module $M$ over a connective $E_1$ ring $\Lambda$ is of \emph{finite Tor-amplitude} if there exists an integer $n$ such that for any right $\Lambda$-module $N$ with $\pi_iN=0$ for $i\neq 0$, one has $\pi_i(N\wedge_{\Lambda}M)=0$ for $i\geq n$.

We will repeatedly use the observation of Lurie \cite[Pr. 8.2.5.23(4)]{HA} that if $\Lambda$ is a left coherent $E_1$ ring, then a left $\Lambda$-module is perfect just in case it is both almost perfect and of finite Tor-amplitude.

In \cite{ku}, Blumberg and Mandell prove for any left coherent $E_1$ ring $\Lambda$, the \emph{$G$-theory} of $\Lambda$ (i.e., the algebraic $K$-theory of the category of coherent left $\Lambda$-modules) agrees with the $G$-theory of $\pi_0\Lambda$. In general, coherent $\Lambda$-modules need not be perfect, but when this does happen, $\Lambda$ is said to be \emph{almost regular}:
\begin{dfn} A left coherent $E_1$ ring $\Lambda$ is said to be \emph{almost regular} if any coherent left $\Lambda$-module has finite Tor-amplitude.
\end{dfn}

To prove that an $E_1$ ring is almost regular, it's enough to check that any coherent left $\Lambda$-module is of finite Tor-amplitude.

\begin{prp} A left coherent $E_1$ ring $\Lambda$ such that $\pi_0\Lambda$ is a regular (commutative) ring is almost regular just in case the left $\Lambda$-module $H\pi_0\Lambda$ is of finite Tor-amplitude (and hence perfect).\label{prp:toramp}
\begin{proof} One implication is clear. Since the class of perfect  modules forms a thick subcategory of left $\Lambda$-modules, it follows that $\Lambda$ is almost regular just in case every almost perfect $\Lambda$-module concentrated in degree $0$ is perfect. In other words, it suffices to show that for any left $\pi_0\Lambda$-module $M$ of finite presentation, the left $\Lambda$-module $HM$ is perfect. Now since $\pi_0\Lambda$ is regular and coherent, $HM$ is perfect when regarded as an $H\pi_0\Lambda$-module, and since the $E_1$ map $\fromto{\Lambda}{H\pi_0\Lambda}$ has finite Tor-amplitude, $HM$ is perfect when regarded as a $\Lambda$-module as well.
\end{proof}
\end{prp}

\begin{dfn} We'll say that a left coherent $E_1$ ring is \emph{regular} if $\pi_0\Lambda$ is a regular (commutative) ring, and $H\pi_0\Lambda$ has finite Tor-amplitude as a $\Lambda$-module.
\end{dfn}

The previous proposition tells us that a regular $E_1$ ring is almost regular. The following proposition allows us to reduce the question of regularity even further to a question about the left $H\pi_0\Lambda$-module $H\pi_0\Lambda\wedge_{\Lambda}H\pi_0\Lambda$.

\begin{prp} Suppose $\Lambda$ a left coherent $E_1$ ring such that $\pi_0\Lambda$ is regular (commutative). Then $\Lambda$ is regular just in case the left $H\pi_0\Lambda$-module $H\pi_0\Lambda\wedge_{\Lambda}H\pi_0\Lambda$ is of finite Tor-amplitude (hence perfect).
\begin{proof} The forward implication is clear.

For the reverse implication, we observe that the left $\Lambda$-module $H\pi_0\Lambda$ is of finite Tor-amplitude, since for any $\Lambda$-module $N$ such that $\pi_iN=0$ for $i\neq 0$, we have
\begin{equation*}
N\wedge_{\Lambda}H\pi_0\Lambda\simeq N\wedge_{H\pi_0\Lambda}(H\pi_0 \Lambda \wedge_\Lambda H\pi_0\Lambda),
\end{equation*}
and since $H\pi_0 \Lambda \wedge_\Lambda H\pi_0\Lambda$ is of finite Tor-amplitude, the proof is complete.
\end{proof}
\end{prp}


\section{D\'evissage and localization} The main value of the almost regularity hypothesis is that it gives a natural map
\begin{equation*}
\fromto{K(\pi_0\Lambda)\simeq G(\pi_0\Lambda)\simeq G(\Lambda)}{K(\Lambda)},
\end{equation*}
since every coherent module is perfect. Note that, for ordinary rings $R$, every perfect module is coherent, and so this form of regularity gives an equivalence $G(R)\simeq K(R)$. The same thing happens for connective $E_1$ rings with only finitely many nontrivial homotopy groups, so a regular $E_1$ ring $\Lambda$ with only finitely many nontrivial homotopy groups has the same $K$-theory as $\pi_0\Lambda$. However, for $E_1$ rings with infinitely many homotopy groups, this isn't the case.

Nevertheless, regularity coupled with the Blumberg--Mandell D\'evissage Theorem gives an identification of $K(\pi_0\Lambda)$ with the $K$-theory of the category of bounded perfect (left) modules on $\Lambda$, whence we obtain the following localization theorem \cite[Th. 8.8]{heart}:
\begin{thm} Suppose $\Lambda$ a regular $E_1$ ring, and suppose
\begin{equation*}
L\colon\fromto{\mathbf{Mod}(\Lambda)}{\mathbf{Mod}(\Lambda)}
\end{equation*}
a smashing localization functor. Assume that a perfect left $\Lambda$-module $M$ is $L$-acyclic just in case it is truncated. Then there is a fiber sequence of infinite loop spaces
\begin{equation}\label{eq:Llocalize}
K(\pi_0\Lambda)\to K(\Lambda)\to K(L(\Lambda)).
\end{equation}
\end{thm}

\begin{nul} Note that if $\Lambda$ is a regular $E_1$ ring and if $L\colon\fromto{\mathrm{Mod}(\Lambda)}{\mathrm{Mod}(\Lambda)}$ is a smashing localization such that a left $\Lambda$-module is $L$-acyclic just in case it is truncated, then $L$ is necessarily finite. Indeed, every truncated $\Lambda$-module can be written as a filtered colimit of coherent $\Lambda$-modules, which are compact objects of $\mathrm{Mod}(\Lambda)$.
\end{nul}

There is a relevant special case to consider here. First, recall that if $\Lambda$ is an $E_1$ ring and if $S\subset\pi_{\ast}\Lambda$ is a multiplicatively closed set of homogeneous elements, then a left $\Lambda$-module $M$ is said to be \emph{$S$-nilpotent} if for every $x\in\pi_{\ast}M$, there exists $s\in S$ such that $sx=0$.

\begin{cor}\label{thm:BM} Suppose $\Lambda$ a regular $E_1$ ring, and suppose $S\subset\pi_{\ast}\Lambda$ a multiplicatively closed collection of homogeneous elements satisfying the left Ore condition \cite[Df. 8.2.4.1]{HA} with the additional property that a left $\Lambda$-module $M$ is $S$-nilpotent just in case it is truncated. Then there is a fiber sequence of infinite loop spaces
\begin{equation}\label{eq:Slocalize}
K(\pi_0\Lambda)\to K(\Lambda)\to K(S^{-1}\Lambda).
\end{equation}
\end{cor}

\begin{nul}\label{nul:localizeaway} More generally, we may consider the following circumstance. If $\Lambda$ is an $E_1$ ring and if $\mathscr{S}=\{S_j\}_{j\in J}$ is a finite family of multiplicatively closed subsets of homogeneous elements of $\pi_{\ast}\Lambda$, then a left $\Lambda$-module $M$ will be said to be $\mathscr{S}$-nilpotent if it is $S_j$-nilpotent for every $j\in J$.

For any nonempty subset $I\subset J$, we can consider the smallest multiplicatively closed set $S_I\subset\pi_{\ast}\Lambda$ of homogenous elements containing $S_i$ for every $i\in I$, and we have the smashing localization functor $L_{S_I}\colon\goesto{M}{S_I^{-1}M}$ on the $\infty$-category $\mathrm{Mod}(\Lambda)$.

Now if $\PP^{\ast}(J)$ is the poset of nonempty subsets of $J$, then the assignment $\goesto{I}{L_{S_I}}$ defines a diagram
\begin{equation*}
\fromto{N\PP^{\ast}(J)}{\mathrm{Fun}(\mathrm{Mod}(\Lambda),\mathrm{Mod}(\Lambda))}
\end{equation*}
of localization functors. The (homotopy) limit of this functor will be denoted $L_{\mathscr{S}}\colon\fromto{\mathrm{Mod}(\Lambda)}{\mathrm{Mod}(\Lambda)}$. Using the characterization of \cite[Pr. 5.2.7.4(3)]{HTT} and the fact that both localization functors and colimits in stable $\infty$-categories are exact, we deduce that $L_{\mathscr{S}}$ is again a smashing localization functor. Furthermore, the $\infty$-category of $L_{\mathscr{S}}$-acyclics is the $\infty$-category of $\mathscr{S}$-nilpotent $\Lambda$-modules. 

When $\Lambda$ is $E_{\infty}$ and each $S_j$ is generated by a homogeneous element $\beta_j\in\pi_{\ast}\Lambda$, the localization $L_{\mathscr{S}}$ agrees with the localization away from the finitely generated ideal $(\beta_j)_{j\in J}$ in the sense of \cite[\S 5]{gmcompletions}.
\end{nul}

\begin{cor}\label{thm:betterBM} Suppose $\Lambda$ a regular $E_1$ ring, and suppose $\mathscr{S}=\{S_j\}_{j\in J}$ a finite family of multiplicatively closed subsets of homogeneous elements of $\pi_{\ast}\Lambda$, and suppose that for each nonempty subset $I\subset J$, the multiplicative subset $S_I\subset\pi_{\ast}\Lambda$ satisfies the left Ore condition \cite[Df. 8.2.4.1]{HA}. Suppose additionally that a left $\Lambda$-module $M$ is $\mathscr{S}$-nilpotent just in case it is truncated. Then there is a fiber sequence of infinite loop spaces
\begin{equation}\label{eq:Jlocalize}
K(\pi_0\Lambda)\to K(\Lambda)\to K(L_{\mathscr{S}}\Lambda).
\end{equation}
\end{cor}

When the first author wrote \cite{heart}, he was unable to verify the almost regularity of various interesting $E_1$ rings, and so was unable to prove that \eqref{eq:Llocalize} was a fiber sequence in many interesting cases. In this paper, armed with the regularity criteria of the first section, we produce such localization sequences for $\KO$, $\TMF$, and $\BP{n}^{\ast}$.


\section{A localization sequence for $\mathrm{KO}$}

\begin{prp} The connective real $K$-theory spectrum $\ko$ is regular.
\begin{proof}
The theorem of Reg Wood \cite[p. 206]{Adams} says that there is a cofiber sequence
\begin{equation*}
\Sigma \ko\ \tikz[baseline]\draw[>=stealth,->,font=\scriptsize](0,0.5ex)--node[above]{$\eta$}(0.5,0.5ex);\  \ko \to \ku
\end{equation*}
of $\ko$-modules. The long exact sequence on homotopy groups implies that there is a cofiber sequence
\begin{equation*}
\Sigma^2 \ku\ \tikz[baseline]\draw[>=stealth,->,font=\scriptsize](0,0.5ex)--node[above]{$\beta$}(0.5,0.5ex);\  \ku \to H\ZZ,
\end{equation*}
where $\beta \in \pi_2(ku)$ is the Bott element.  Therefore, $H\ZZ$ is a perfect $\ko$-module.
\end{proof}
\end{prp}

Now since a perfect $\ko$-module $M$ is nilpotent with respect to the periodicity generator $\beta^4 \in \pi_8(ko)$ just in case it is truncated, the previous result gives the following special case of \eqref{eq:Slocalize}.
\begin{cor} The canonical $E_{\infty}$ ring morphisms $\fromto{\ko}{H\ZZ}$ and $\fromto{\ko}{\KO}$ give a $K$-theory fiber sequence
\begin{equation*}
K(\ZZ)\to K(\ko)\to K(\KO),
\end{equation*}
recovering the result discussed in the introduction.
\end{cor}

\begin{nul} This confirms a conjecture posed by Christian Ausoni and John Rognes in \cite[Rk. 3.7]{rational}. Indeed, as a consequence of their rational computation of $K(\ko)$, we obtain a rational computation of $K(\KO)$. In particular, Borel's computation of $K_{\ast}(\ZZ)\otimes\QQ$ implies that the Poincar\'e series for $K(\ZZ)$ is
\begin{equation*}
1+\frac{t^5}{1-t^4},
\end{equation*}
and in \cite[Rk. 3.4]{rational} it is shown that the Poincar\'e series for $K(\ko)$ is
\begin{equation*}
1+\frac{2t^5}{1-t^4};
\end{equation*}
hence the Poincar\'e series for $K(\KO)$ is
\begin{equation*}
1+t+\frac{2t^5+t^6}{1-t^4},
\end{equation*}
which agrees with the expectation of Ausoni--Rognes.

Furthermore, it is observed in loc. cit. that the natural map
\begin{equation*}
\fromto{K(\ko)\wedge H\QQ}{(K(\ku)\wedge H\QQ)^{h\ZZ/2}}
\end{equation*}
is a rational equivalence. Since the action of $\ZZ/2$ on $K(\ZZ)$ is trivial, we obtain a morphism of fiber sequences
\begin{equation*}
\begin{tikzpicture} 
\matrix(m)[matrix of math nodes, 
row sep=4ex, column sep=4ex, 
text height=1.5ex, text depth=0.25ex] 
{K(\ZZ)\wedge H\QQ&K(\ko)\wedge H\QQ&K(\KO)\wedge H\QQ\\ 
F(B(\ZZ/2)_+,K(\ZZ)\wedge H\QQ)&(K(\ku)\wedge H\QQ)^{h\ZZ/2}&(K(\KU)\wedge H\QQ)^{h\ZZ/2}\\}; 
\path[>=stealth,->,font=\scriptsize] 
(m-1-1) edge (m-1-2) 
edge (m-2-1) 
(m-1-2) edge (m-1-3)
edge (m-2-2)
(m-1-3) edge (m-2-3)
(m-2-1) edge (m-2-2)
(m-2-2) edge (m-2-3); 
\end{tikzpicture}
\end{equation*}
where the lower fiber sequence is the homotopy fixed-point sequence of the rationalization of Blumberg-Mandell's localization sequence \cite{ku}. Since the left-hand vertical map is an equivalence, it follows that, similarly, the natural map
\begin{equation*}
\fromto{K(\KO)\wedge H\QQ}{(K(\KU)\wedge H\QQ)^{h\ZZ/2}}
\end{equation*}
is a rational equivalence.  This confirms the expectation of Ausoni--Rognes, and confirms in this case the suspicion that rationalized $K$-theory should satisfy Galois descent, even in the derived algebro-geometric context.
\end{nul}


\section{A localization sequence for $\TMF$}

\begin{prp} The topological modular forms spectrum $\tmf$ is regular.
\begin{proof}
We have shown in Proposition~\ref{prp:toramp} that it is necessary and sufficient that $H\ZZ \wedge_\tmf H\ZZ$ be of finite Tor-amplitude. As $\ZZ$ has finite projective dimension, this is the case if and only if it is truncated.

The ring $\pi_{\ast}\tmf[1/6]$ is $\ZZ[1/6,c_4,c_6]$ and $\tmf$ is commutative, so we may construct cofibration sequences of $\tmf$-modules
\[
\Sigma^8 \tmf[1/6]\ \tikz[baseline]\draw[>=stealth,->,font=\scriptsize](0,0.5ex)--node[above]{$c_4$}(0.5,0.5ex);\ \tmf[1/6] \to \tmf[1/6]/c_4
\]
and
\[
\Sigma^{12} \tmf[1/6]/c_4\ \tikz[baseline]\draw[>=stealth,->,font=\scriptsize](0,0.5ex)--node[above]{$c_6$}(0.5,0.5ex);\ \tmf[1/6]/c_4 \to H\ZZ[1/6].
\]
Therefore, $H\ZZ \wedge_\tmf H\ZZ$ has trivial homotopy groups above degree $22$ after inverting $6$.

The higher homotopy groups therefore consist of $6$-torsion.  Applying the long exact sequence associated to the cofibration sequence $H\ZZ \to H\ZZ \to H\ZZ/p$ twice each for the primes $2$ and $3$, we find that the higher homotopy groups vanish if and only if the higher homotopy groups of $H\ZZ/2 \wedge_\tmf H\ZZ/2$ and $H\ZZ/3 \wedge_\tmf H\ZZ/3$ vanish.

At the prime $3$, Hill \cite{tmfat3} has determined the Hopf algebra
of homotopy groups of $H\ZZ/3 \wedge_\tmf H\ZZ/3$, and in particular
shown that it is finite-dimensional over $\ZZ/3$.

At the prime $2$, Mathew \cite{tmfat2} has shown that the map
\begin{equation*}
H\ZZ/2 \wedge \tmf \to H\ZZ/2 \wedge H\ZZ/2
\end{equation*}
is the inclusion of the subalgebra $\ZZ/2[\xxi_1^8,\xxi_2^4,\xxi_3^2,\xxi_4,\ldots]$ of the dual Steenrod algebra; in particular, the dual Steenrod algebra is a free module over it.  By taking the identity of smash products
\[
H\ZZ/2 \wedge_\tmf H\ZZ/2 \simeq (H\ZZ/2 \wedge H\ZZ/2) \wedge_{H\ZZ/2 \wedge \tmf} H\ZZ/2
\]
and applying the K\"unneth spectral sequence of \cite[IV 4.1]{ekmm},
we find that $\pi_* (H\ZZ/2 \wedge_\tmf H\ZZ/2)$ is the quotient algebra $\ZZ/2[\xxi_1,\xxi_2,\xxi_3]/(\xxi_1^8,\xxi_2^4,\xxi_3^2)$.  In particular, it is finite-dimensional.
\end{proof}
\end{prp}

We note that \cite{tmfat2}, based on work of Hopkins-Mahowald, also establishes that there is an $8$-cell complex $DA(1)$ such that, $2$-locally, $\tmf \wedge DA(1)$ is a form of $\BP{2}$.  However, to make use of this as we did for $ko$ we would need to extend this to a \emph{multiplicative} identification.

In any case, we have the following consequence of \eqref{eq:Jlocalize}.
\begin{cor} The canonical $E_{\infty}$ ring morphisms $\fromto{\tmf}{H\ZZ}$ and $\fromto{\tmf}{\TMF}$ give a $K$-theory fiber sequence
\begin{equation*}
K(\ZZ)\to K(\tmf)\to K(\TMF).
\end{equation*}
\begin{proof}
Let $L=L_{\mathscr{S}}$ (in the notation of \ref{nul:localizeaway}) for the finite family
\begin{equation*}
\mathscr{S}=\{\langle c_4\rangle,\langle \Delta\rangle\},
\end{equation*}
where $\langle x\rangle\subset\pi_{\ast}\tmf$ is the multiplicative subset of homogeneous elements generated by $x$. A perfect $\tmf$-module is truncated just in case it is both $c_4$-nilpotent and $\Delta$-nilpotent. By the previous results, it suffices to show that $\TMF = L_{\mathscr{S}} (\tmf)$.

The map $\tmf \to \TMF$ is a connective cover.  As the elements of $\mathscr{S}$ are generated by elements of positive degree, this map is automatically an $L_{\mathscr{S}}$-equivalence.

We recall that $\TMF$ is the global section object $\Gamma(\overline{\mathcal{M}},\mathscr{O}^{der})$ of a sheaf of $E_\infty$ rings on the compactified moduli of elliptic curves \cite{buildtmf}. The vanishing loci of $c_4$ and $\Delta$ have empty intersection on $\overline{\mathcal{M}}$, and so we obtain a cover
\[
(c_4^{-1}\overline{\mathcal{M}}) \amalg (\Delta^{-1} \overline{\mathcal{M}}) \to \overline{\mathcal{M}}.
\]
Descent for $\mathscr{O}^{der}$ with respect to the associated \v{C}ech cover expresses $\TMF$ as the homotopy pullback of the diagram
\[
c_4^{-1} \TMF \rightarrow (c_4 \Delta)^{-1} \TMF \leftarrow \Delta^{-1} \TMF.
\]
In particular, $\TMF$ is $\mathscr{S}$-local.
\end{proof}
\end{cor}



\section{A localization sequence for $\BP{n}$} Fix a prime $p$ and a
generalized truncated Brown--Peterson spectrum $\BP{n}$ \cite[\S
3]{bp2} with an $E_1$ $MU$-algebra structure \cite[p. 506]{MR1990937}.
As $\pi_* \BP{n}$ has finite projective dimension, a simple argument
using the K\"unneth spectral sequence immediately implies that the
$E_1$ ring $\BP{n}$ is regular.

Let $L=L_{\mathscr{S}}$ (in the notation of \ref{nul:localizeaway}) for the finite family
\begin{equation*}
\mathscr{S}=\{\langle v_1\rangle,\dots,\langle v_n\rangle\},
\end{equation*}
where $\langle v_i\rangle\subset\pi_{\ast}\BP{n}$ is the multiplicative subset of homogeneous elements generated by $v_i$. 

Equivalently, if $J$ denotes the ideal $(v_1,\ldots,v_n)$, and one selects lifts of these elements to $v_i \in MU_*$, then $L$ is the localization away from $J$ in the sense of \cite[\S 5]{gmcompletions}.  Either by the discussion of \ref{nul:localizeaway} or by \cite[5.2]{gmcompletions}, this is a smashing localization on $\BP{n}$-modules. Write $\BP{n}^{\ast}\coloneq L(\BP{n})$.

\begin{prp} A perfect $\BP{n}$-module is $L$-acyclic if and only if it is truncated.
\begin{proof}
By \cite[5.1]{gmcompletions}, a left $\BP{n}$-module is $L$-acyclic if and only if it is $v_k$-torsion for all $1 \leq k \leq n$.  The ring $\pi_* \BP{n}$ is noetherian, so any perfect module $M$ has $\pi_* M$ a finitely generated module over it.  Therefore, $M$ is $L$-acyclic if and only if $\pi_*M$ has a finite filtration whose subquotients are finitely generated $\ZZ_{(p)}$-modules, which occurs precisely when $M$ is truncated.
\end{proof}
\end{prp}

\noindent Note that as a consequence of this characterization, the functor $L$ is actually independent of the choices of the lifts $v_i$.

\begin{cor} The $E_1$ ring morphisms $\fromto{\BP{n}}{H\ZZ_{(p)}}$ and $\fromto{\BP{n}}{\BP{n}^{\ast}}$ give a $K$-theory fiber sequence
\begin{equation*}
K(\ZZ_{(p)})\to K(\BP{n})\to K(\BP{n}^{\ast}).
\end{equation*}
\end{cor}

\bibliography{almost}
\bibliographystyle{amsalpha}

\end{document}